\documentclass[11pt]{amsart}
\usepackage{amsmath}
\usepackage{amscd}

\theoremstyle{definition}
\newtheorem{definition}{Definition}[section]
\newtheorem{defiprop}[definition]{Proposition-Definition}

\theoremstyle{plain}
\newtheorem{theorem}[definition]{Theorem}
\newtheorem{corollary}[definition]{Corollary}
\newtheorem{lemma}[definition]{Lemma}
\newtheorem{proposition}[definition]{Proposition}

\newtheorem*{theorem*}{Theorem}
\newtheorem*{lemma*}{Lemma}
\newtheorem*{proposition*}{Proposition}
\newtheorem*{problem}{Problem}

\newtheorem*{definition*}{Definition}

\numberwithin{equation}{section}

\newcommand{\Proof}{\noindent {\it Proof.\ }}
\newcommand{\dotdot}{{\, . \, . \,}}

\newcommand{\F}{{\mathcal F \!}}

\newcommand{\R}{{\mathbb R}}

\newcommand{\vol}{\mbox{$\mathrm{vol}$}}

\newcommand{\diag}{\mbox{$\mathrm{diag}$}}

\newcommand{\innerxy}[2]{{\langle #1 , #2 \rangle}}
\newcommand{\nmx}[1]{\mbox{$\|#1\|$}}

\newcommand{\ba}{{\bf a}}
\newcommand{\bx}{{\bf x}}
\newcommand{\by}{{\bf y}}

\newcommand{\bv}{{\bf v}}

\newcommand{\be}{{\bf e}}

\newcommand{\bzero}{{\bf 0}}
\newcommand{\ssM}{{\scriptscriptstyle M}}
\newcommand{\ssN}{{\scriptscriptstyle N}}
\newcommand{\bracket}[1] {{<\! #1 \!>}}
\newcommand{\kvector}[2]
    {{\bf #1}_{1} \wedge {\bf #1}_{2} \wedge \ldots \wedge {\bf #1}_{#2}}
\newcommand{\geod}[1]{{\rm Geod}(#1)}

\title{What is wrong with the Hausdorff measure in Finsler spaces}

\author{J.C. \'Alvarez Paiva}
\address{J.C. \'Alvarez Paiva, Department of Mathematics,
         Polytechnic University, Six Metro\-Tech Center, Brooklyn,
         New York, 11201, U.S.A.}
\email{jalvarez@duke.poly.edu}
\thanks{The first author was partially funded by FAPESP grant ${\rm N}^{\rm o}$
         2004/01509-0}

\author{G. Berck}
\address{G. Berck, Universit\'e Catholique de Louvain,
Institut de Math\'ematique Pure et Appl., Chemin du Cyclotron 2, B--1348
Louvain--la--Neuve, Belgium.}
\email{g.berck@math.ucl.ac.be}

\keywords{Finsler manifold, Hausdorff measure, minimal submanifolds,
integral geometry, Holmes-Thompson volume.}

\subjclass{53B40; 49Q05, 53C55}

\begin{document}

\begin{abstract}
We construct a class of  Finsler metrics in three-dimensional space such
that all their geodesics are lines, but not all planes are extremal for
their Hausdorff area functionals. This shows that if the Hausdorff measure
is used as notion of volume on Finsler spaces, then totally geodesic submanifolds
are not necessarily minimal, filling results such as those of
Ivanov~\cite{Ivanov:2002} do not hold, and integral-geometric formulas do not
exist. On the other hand, using the Holmes-Thompson definition of volume, we
prove a general Crofton formula for Finsler spaces and give an easy proof that
their totally geodesic hypersurfaces are minimal.
\end{abstract}

\maketitle

\section{Introduction}

One of the first questions that arise in the study of Finsler manifolds is
whether it is possible to define their volume in a natural way. In
\cite{Busemann:1950} Busemann argues strongly that the volume of a Finsler
manifold should be its Hausdorff measure. His argument is based on a
number of axioms that any natural definition of volume on Finsler
spaces must satisfy. Alas, there are other natural definitions of volume that
satisfy those axioms. Among them is one that is rapidly becoming
the uncontested definition of volume for Finsler spaces: the Holmes-Thompson
volume.

\begin{definition}
The {\it Holmes-Thompson volume\/} of an $n$-dimensional Finsler manifolds equals
the symplectic volume of its unit co-disc bundle divided by the volume of
the Euclidean unit ball of dimension $n$. The {\it area of a submanifold\/} is the
Holmes-Thompson volume for the induced Finsler metric.
\end{definition}

Its deep ties to convexity, differential geometry, and integral geometry
(see \cite{Thompson:1996}, \cite{Alvarez-Thompson:survey}, and
\cite{Alvarez:survey}) have made the Holmes-Thompson volume easier to study
than the Hausdorff measure. However, since many of the theorems that are known
to hold for the Holmes-Thompson volume have not been proved nor disproved for
the Hausdorff measure, one may wonder whether it is just
not suitable for the study of Finsler spaces, or whether the right techniques
for its study have not yet been found. Our results suggest that the Hausdorff
measure is not suitable.

In this paper we show through an experiment that totally geodesic submanifolds
of a Finsler space do not necessarily extremize the Hausdorff area integrand.
The experiment also shows that the Holmes-Thompson volume cannot be replaced
by the Hausdorff measure in the filling result of S. Ivanov \cite{Ivanov:2002},
and that there is no Crofton formula for the Hausdorff measure of hypersurfaces
in Finsler spaces. In \cite{Schneider:2001a} Schneider gives examples of normed
spaces for which there is no Crofton formula for the Hausdorff measure of
hypersurfaces, but his examples are not Finsler spaces---their unit spheres are
not smooth hypersurfaces with positive Gaussian curvature. The phenomenon we
uncover is completely different: the obstruction to the existence of a Crofton
formula does not lie in the normed spaces that make up the tangent bundle, but
in the way they combine to make a Finsler metric. In fact, we can find Finsler
metrics in $\R^3$ whose restriction to a large ball is arbitrarily close to the
Euclidean metric and for which there is no Crofton formula for the Hausdorff
measure of hypersurfaces contained in the ball.

Using the standard norm and inner product in $\R^3$ and identifying
$\R^3 \times \R^3$ with the tangent bundle of $\R^3$, the experiment
can be summarized in the following theorem.

\begin{theorem*}
For any real number $\lambda$, all geodesics of the Finsler metric
$$
\varphi_\lambda (\bx,\bv) =
\frac{(1 + \lambda^2 \nmx{\bx}^2)\nmx{\bv}^2 + \lambda^2 \innerxy{\bx}{\bv}^2}
{\nmx{\bv}}
$$
are straight lines. However, the only value of $\lambda$ for which
all planes are extremals of the Hausdorff  area functional of $\varphi_\lambda$
is $\lambda = 0$.
\end{theorem*}

There is a lot more in this theorem that meets the eye. For example, it takes
some work in Section~\ref{Finsler-section} just to see that the $\varphi_\lambda$
are Finsler metrics. The computation of area integrands for Finsler metrics in
$\R^3$ requires the computation of a Funk transform. It is only due to the
simple formula for $\varphi_\lambda$ that in Section~\ref{integrands-section}
we are able to compute its area integrand in $\R^3 \times \Lambda^2 \R^3$:
$$
\phi_{\lambda}(\bx,\ba) =
\frac{2(1 + \lambda^2\nmx{\bx}^2)^{3/2}
  \left((1 + 2\lambda^2\nmx{\bx}^2)\nmx{\ba}^2 -
                         \lambda^2 \innerxy{\bx}{\ba}^2 \right)^{3/2}}
 {(2 + 3 \lambda^2\nmx{\bx}^2)\nmx{\ba}^2 -  \lambda^2 \innerxy{\bx}{\ba}^2} .
$$
Apparently, this is the first example of an explicit computation of the
(Hausdorff) area integrand of a non-Riemannian Finsler metric.

In Section~\ref{Hamel-section} we characterize smooth parametric integrands
of degree $n-1$ on $\R^n$ for which hyperplanes are extremal as those that
satisfy a certain linear differential equation. The simplicity of this
equation makes it possible to show that all planes cannot be extremals for
the Hausdorff area functional of $\varphi_\lambda$ unless $\lambda$ is zero.
Also in Section~\ref{Hamel-section} we quickly show that Ivanov's filling theorem
implies that two-dimensional totally geodesic submanifolds of Finsler spaces
are minimal for the Holmes-Thompson area functional, and conclude that filling
theorems of this nature cannot hold for the Hausdorff measure.

The proof of the Crofton formula for hypersurfaces in Finsler spaces, its
application to the minimality of totally geodesic hypersurfaces for the
Holmes-Thompson volume, and the proof that there is no Crofton formula
for the Hausdorff area functional of $\varphi_\lambda$ unless $\lambda$ is
zero are all in Section~\ref{Crofton-section}

\section{A class of Finsler metrics}\label{Finsler-section}

Roughly speaking, a Finsler manifold is a manifold in which each tangent space
has been provided with a norm and these norms vary smoothly with the base point.
In order to guarantee that the subject stays in the realm of differential
geometry, it is standard (see \cite{Bao-Chern-Shen}) to ask that the norms be
Minkowski norms.

\begin{definition}
Let $V$ be a finite-dimensional vector space over the reals. A norm
$\varphi$ on $V$ is said to be a {\it Minkowski norm\/} if it is smooth
away from the origin and the Hessian of $\varphi^2$ is positive definite at
every nonzero point.
\end{definition}

Two useful properties of these norms are that they are closed under linear
transformations, and that a norm is a Minkowski norm if and only if its
restrictions to two-dimensional subspaces are Minkowski norms. The following
theorem introduces the class of Minkowski norms that will be used in this paper.

\begin{theorem}\label{Minkowski-norms}
Let $A$ and $B$ be two $n \times n$ symmetric, positive definite matrices.
If $\lambda_1 \leq \cdots \leq \lambda_n$ are the roots of
$\det(A - \lambda B) = 0$---the eigenvalues of $A$ relative to $B$---,
then $F(\bv) = \innerxy{A\bv}{\bv}/\sqrt{\innerxy{B\bv}{\bv}}$ is a Minkowski norm
on $\R^n$ if and only if $\lambda_n < 2\lambda_1$.
\end{theorem}

\Proof
We first prove the theorem when $n=2$. By a linear change of variables, we
may assume that
$$
F(v_1,v_2) = \frac{\lambda_1 v_1^2 + \lambda_2 v_2^2}{\sqrt{v_1^2 + v_2^2}} .
$$
Clearly $F$ is nonnegative and homogeneous of degree one. Because we are working
in two dimensions, we can  verify that the Hessian of $F^2$ is positive definite
at every nonzero point simply by checking that the Laplacian of $F$ is positive
away from the origin. A computation shows that
$$
\Delta F(v_1,v_2) := \frac{(2\lambda_2 - \lambda_1)v_1^2 +
                     (2\lambda_1 - \lambda_2)v_2^2}{(v_1^2 + v_2^2)^{3/2}} ,
$$
which is positive away from the origin if and only if $\lambda_2 < 2\lambda_1$.
This settles the two-dimensional case.

Now we turn to the general case. After a suitable linear change of variable, we
may assume that $B$ is the identity matrix and that $A$ is the diagonal matrix
$\diag(\lambda_1,\ldots,\lambda_n)$.

By the Rayleigh-Courant-Fisher theorem (see \cite[pp. 110--113]{Arnold:1989})
the eigenvalues $\lambda_1' \leq \lambda_2'$ of the restriction of the quadratic
form $\lambda_1 v_1^2 + \cdots \lambda_n v_n^2$ to a two-dimensional subspace
satisfy the inequalities $\lambda_1 \leq \lambda_1'$ and
$\lambda_2' \leq \lambda_n$. Thus, if $\lambda_n \leq 2\lambda_1$, then
$\lambda_2' < 2\lambda_1'$ and the restriction of $F$ to every two-dimensional
subspace is a Minkowski norm. Therefore, $F$ is a Minkowski norm on $\R^n$.

On the other hand if $F$ is a Minkowski norm, its restriction to the
two-dimensional subspace spanned by the basis vectors $\be_1$ and
$\be_n$ is also a Minkowski norm and, therefore, $\lambda_n < 2 \lambda_1$.
\qed \medskip

A {\it Finsler metric\/} on a manifold is simply a continuous function
on its cotangent bundle that is smooth away from the zero section and such
that its restriction to each tangent space is a Minkowski norm. The special
class of Finsler metrics we will be working with are of the following form:

\begin{theorem}\label{Finsler-metrics}
Let $g$ and $h$ be two Riemannian metrics on a manifold $M$. The function
$\varphi := g/\sqrt{h}$ is a Finsler metric if and only if on every tangent
unit sphere of $(M,h)$ the maximum of the function $\bv \mapsto g(\bv,\bv)$
is less than twice its minimum.
\end{theorem}

\Proof
The function $\varphi$ is clearly continuous and smooth away from the zero
section. It remains for us to see that its restriction to an arbitrary tangent
space $T_{\bx}M$ is a Minkowski norm. For this, consider the quadratic forms
$g_{\bx}$ and $h_{\bx}$ on $T_{\bx}M$. The variational characterization of
eigenvalues of positive-definite quadratic forms implies that the smallest and
the largest of the eigenvalues of $g_\bx$ relative to $h_\bx$ are, respectively,
the minimum and the maximum of the restriction of $g_\bx$ to the hypersurface
$h_\bx = 1$. Theorem~\ref{Minkowski-norms} tells us that $g_\bx /\sqrt{h_\bx}$
is a Minkowski norm if and only if this maximum is less than twice the minimum.
\qed

\begin{corollary}\label{Main-examples}
For every value of the parameter $\lambda$, the function
$$
\varphi_\lambda (\bx,\bv) =
\frac{(1 + \lambda^2 \nmx{\bx}^2)\nmx{\bv}^2 + \lambda^2 \innerxy{\bx}{\bv}^2}
{\nmx{\bv}}
$$
is a Finsler metric on $\R^n$.
\end{corollary}

\Proof
By Theorem~\ref{Finsler-metrics}, it suffices to verify that for each fixed
$\bx$ the maximum of the restriction of
$(1 + \lambda^2 \nmx{\bx}^2)\nmx{\bv}^2 + \lambda^2 \innerxy{\bx}{\bv}^2$
to the sphere $\nmx{\bv} = 1$ is less than twice its minimum. This is indeed
the case:
$$
(1 + \lambda^2 \nmx{\bx}^2) + \lambda^2 \nmx{\bx}^2 < 2(1 + \lambda^2 \nmx{\bx}^2) .
$$
\qed

\section{Hausdorff area integrands}\label{integrands-section}

The computation of (Hausdorff) area integrands on Finsler manifolds depends on the
following key result of Busemann \cite{Busemann:1947}.

\begin{theorem}
The Hausdorff measure on an $n$-dimensional normed space is characterized
by its invariance under translations and by the fact that the measure of
the unit ball equals the volume of the unit ball in an $n$-dimensional
Euclidean space.
\end{theorem}

The theorem implies that if $V$ is an $n$-dimensional normed space with
unit ball $B$, the $k$-area integrand $(1 \leq k \leq n)$ of $V$ is the function
that assigns to every nonzero simple $k$-vector $a = \kvector{v}{k}$ the number
$\epsilon_{k}\vol(B \cap \bracket{a};a)^{-1}$. Here $\epsilon_k$ denotes the
volume of the Euclidean unit ball of dimension $k$, $\bracket{a}$ denotes the
subspace spanned by the vectors $\bv_1,\ldots,\bv_k$, and
$\vol(B \cap \bracket{a};a)$ denotes the $k$-dimensional volume of the
intersection of $B$ and the subspace $\bracket{a}$ with respect to the Lebesgue
measure on $\bracket{a}$ determined by the basis $\bv_{1},\ldots,\bv_{k}$.
Busemann showed in \cite{Busemann:1949a} that the $(n-1)$-area integrand
of an $n$-dimensional normed space $V$ is a norm on $\Lambda^{(n-1)} V$.

To avoid confusion with the Holmes-Thompson area integrands, we shall use
the term {\it Hausdorff $k$-area integrands\/} for the integrands described
in the preceding paragraph. In \cite{Alvarez-Thompson:survey} the Hausdorff
$k$-area integrands are called {\it Busemann $k$-volume densities.}

Since our examples are three-dimensional, we will concentrate on the effective
computation of the Hausdorff $2$-area integrand of a norm on $\R^3$. It will
be useful to consider the natural Euclidean structures on $\R^3$ and
$\Lambda^2 \R^3$. On this last space the standard (orthonormal) basis is
formed by the vectors $\be_2 \wedge \be_3$, $\be_3 \wedge \be_1$, and
$\be_1 \wedge \be_2$, where $\be_1$, $\be_2$, and $\be_3$ form the standard
basis of $\R^3$. The coordinates $(a_1,a_2,a_3)$ will represent the bivector
$\ba = a_1 \be_2 \wedge \be_3 + a_2 \be_3 \wedge \be_1 + a_3 \be_1 \wedge \be_2$,
and will be used to identify $\Lambda^2 \R^3$ with $\R^3$. With this
identification, and identifying $\Lambda^3 \R^3$ with $\R$ via the basis
$\be_1 \wedge \be_2 \wedge \be_3$, a three-vector $\bv \wedge \ba$ is
identified with the scalar $\innerxy{\bv}{\ba}$.

Two concepts that are useful in the computation of Hausdorff area integrands
are the radial function of a norm on a Euclidean space and the Funk transform of
a function on the Euclidean $2$-sphere.

\begin{definition}
The {\it radial function\/} of a norm $\varphi$ on a Euclidean space is
the restriction of $\varphi^{-1}$ to the Euclidean unit sphere.
\end{definition}

Note that the norm can be easily recovered from the radial function through
the formula $\varphi(\bx) = \nmx{\bx}/\rho(\bx/\nmx{\bx})$.

\begin{definition}
If $f$ is a continuous real-valued function on $S^2$, its {\it Funk transform\/}
is the function on the sphere given by
$$
\F f(\bx) = \int_{\innerxy{\bx}{\bv} = 0} \! f \, d\sigma ,
$$
where $\innerxy{\bx}{\bv} = 0$ denotes the great circle perpendicular to
$\bx$ and $d\sigma$ is its element of arclength.
\end{definition}

\begin{proposition}\label{Funk-transform}
If $\rho : S^2 \rightarrow (0 \dotdot \infty)$ is the radial function of
a norm on $\R^3$, the radial function of its Hausdorff $2$-area integrand is
the Funk transform of $\rho^2$ divided by $2 \pi$.
\end{proposition}

\Proof
If $B$ is the unit ball of a norm on $\R^3$ and $\ba$ is a unit bivector
for the Euclidean norm on $\Lambda^2 \R^3$, the quantity
$\vol(B \cap \bracket{\ba}; \ba)$ is just the Euclidean area
of the intersection of $B$ and the plane
$\{\bv \in \R^3 : \ba \wedge \bv = \bzero \}$. Using polar coordinates
$$
\vol(B \cap \bracket{\ba}; \ba) =
\frac{1}{2} \int_{\ba \wedge \bv = \bzero} \rho^2 \, d\sigma =
\F \rho^2 (\ba)/2 .
$$
Since the value of the Hausdorff $2$-area integrand on $\ba$ is
$\pi/\vol(B \cap \bracket{\ba};\ba)$, the radial function is
$\ba \mapsto \F \rho^2 (\ba)/2\pi$.
\qed \medskip

In some cases, the previous proposition allows the explicit computation of
Hausdorff area integrands.

\begin{theorem}\label{area-integrand}
The Hausdorff $2$-area integrand of the Finsler metric
$$
\varphi_\lambda (\bx,\bv) =
\frac{(1 + \lambda^2 \nmx{\bx}^2)\nmx{\bv}^2 + \lambda^2 \innerxy{\bx}{\bv}^2}
{\nmx{\bv}}
$$
on $\R^3$ is
$$
\phi_{\lambda}(\bx,\ba) =
\frac{2(1 + \lambda^2\nmx{\bx}^2)^{3/2}
  \left((1 + 2\lambda^2\nmx{\bx}^2)\nmx{\ba}^2 -
                         \lambda^2 \innerxy{\bx}{\ba}^2 \right)^{3/2}}
 {(2 + 3 \lambda^2\nmx{\bx}^2)\nmx{\ba}^2 -  \lambda^2 \innerxy{\bx}{\ba}^2} .
$$
\end{theorem}

The main computation is carried out in the following lemma:

\begin{lemma}
Let $\lambda$ be a real number, $s$ a positive real number, and $\bx$ a vector
in $\R^3$. The Funk transform of the function
$f(\bv) = (s^2 + \lambda^2 \innerxy{\bx}{\bv}^2)^{-2}$ is the function
$$
\F f(\ba) = \frac{\pi}{s^3}
\frac{2s^2 + \lambda^2 \|\bx\|^2 - \lambda^2 \innerxy{\bx}{\ba}^2}
     {(s^2 + \lambda^2 \|\bx\|^2 - \lambda^2 \innerxy{\bx}{\ba}^2)^{3/2}} .
$$
\end{lemma}

\Proof
We wish to compute the integral
\begin{equation}\label{int-one}
\F f(\ba) = \int_{\innerxy{\ba}{\bv} = 0}
\frac{d\sigma}{(s^2 + \lambda^2 \innerxy{\bx}{\bv}^2)^{2}}
\end{equation}
for a fixed unit (bi)vector $\ba$. If $\by = \bx - \innerxy{\bx}{\ba}\ba$
is the projection of $\bx$ onto the plane orthogonal to $\ba$, and $\bv$
is a unit vector on that plane, then
$$
\innerxy{\bx}{\bv} = \innerxy{\by}{\bv} = \nmx{\by}\cos(\theta) ,
$$
where $\theta$ is the angle formed by $\by$ and $\bv$. The integral
in~(\ref{int-one}) may now be written as
$$
\F f(\ba) = \int_{0}^{2\pi}
\frac{d\theta}{(s^2 + \lambda^2 \nmx{\by}^2 \cos(\theta)^2)^{2}}
$$
or, upon setting $\kappa = \lambda \nmx{\by}/s$, as
\begin{equation}\label{int-two}
\F f(\ba) = \frac{1}{s^4} \int_{0}^{2\pi}
\frac{d\theta}{(1 + \kappa^2 \cos(\theta)^2)^{2}} .
\end{equation}
Note that the integrand in~(\ref{int-two}) is the pull-back of the closed
form
$$
\chi := \frac{(x^2 + y^2) (x \, dy - y \, dx)}
{\left((1 + \kappa^2)x^2 + y^2 \right)^2}
$$
under the map $\theta \mapsto (\cos(\theta),\sin(\theta))$. Taking $\gamma$ to
be the parameterized ellipse
$\theta \mapsto \left((1 + \kappa^2)^{-1/2}\cos(\theta),\sin(\theta)\right)$
and using Stokes' theorem, it is easy to compute that
$$
\int_{0}^{2\pi} \frac{d\theta}{(1 + \kappa^2 \cos(\theta)^2)^{2}} =
\int_\gamma \chi = \frac{\pi (2 + \kappa^2)}{(1 + \kappa^2)^{3/2}} .
$$
This pretty way to evaluate the integral in~(\ref{int-two}) was suggested to
the authors by Thomas P\"uttmann.

Substituting back $\kappa = \lambda \nmx{\by}/s$ and
$\nmx{\by}^2 = \nmx{\bx}^2 - \innerxy{\bx}{\ba}$ yields
\begin{eqnarray*}
\F f(\ba) &=& \frac{1}{s^4} \int_{0}^{2\pi}
\frac{d\theta}{(1 + \kappa^2 \cos(\theta)^2)^{2}} \\
&=&
\frac{\pi}{s^3}
\frac{2s^2 + \lambda^2 \|\bx\|^2 - \lambda^2 \innerxy{\bx}{\ba}^2}
     {(s^2 + \lambda^2 \|\bx\|^2 - \lambda^2 \innerxy{\bx}{\ba}^2)^{3/2}} .
\end{eqnarray*}
\qed \medskip

\noindent {\it Proof of Theorem~\ref{area-integrand}.\/}
Fixing $\bx$, we see that the radial function of the restriction of
$\varphi_\lambda$ to $T_\bx \R^3$ is
$f(\bv) = (1 + \lambda^2 \nmx{\bx}^2 + \lambda^2 \innerxy{\bx}{\bv}^2)^{-1}$.
By Proposition~\ref{Funk-transform}, the radial function of the Hausdorff
$2$-area integrand in $\Lambda^2 T_\bx \R^3$ is $\rho := \F f^2 /2\pi$. Applying
the preceding lemma with $s = 1 + \lambda^2 \nmx{\bx}^2$, we see that
\begin{equation}
\rho(\ba) =
\frac{2 + 3\lambda^2 \nmx{x}^2 - \lambda^2 \innerxy{\bx}{\ba}^2}
{2(1 + \lambda^2 \nmx{\bx}^2)^{3/2}
  (1 + 2\lambda^2 \nmx{\bx}^2 - \lambda^2 \innerxy{\bx}{\ba}^2)^{3/2}} .
\end{equation}
It follows that if  $\ba$ is any nonzero bivector,
\begin{eqnarray*}
\phi_\lambda (\ba) &=& \nmx{\ba}/\rho(\ba/\nmx{\ba}) \\
&=&
\frac{2(1 + \lambda^2\nmx{\bx}^2)^{3/2}
  \left((1 + 2\lambda^2\nmx{\bx}^2)\nmx{\ba}^2 -
                         \lambda^2 \innerxy{\bx}{\ba}^2 \right)^{3/2}}
 {(2 + 3 \lambda^2\nmx{\bx}^2)\nmx{\ba}^2 -  \lambda^2 \innerxy{\bx}{\ba}^2} .
\end{eqnarray*}

\qed

\section{Generalized Hamel's equations}\label{Hamel-section}

Following Federer \cite{Federer:1969} we define a parametric integrand of degree
$k$ $(1 \leq k \leq n)$ on $\R^n$ as a continuous function that is positively
homogenous of degree one in its second variable
($\Phi(\bx,ta) = t\Phi(\bx,a)$ for $t > 0$).
By analogy with the differential $k$-forms they generalize, parametric
integrands of degree $k$ that are smooth away from
$\R^n \times \{\bzero\}$ will be called {\it differential $k$-integrands.}

A differential $k$-integrand $\Phi$ is said to be {\it projective\/} if
$k$-dimensional flats are extremals of the functional
$$
N \longmapsto \int_N \Phi .
$$
These integrands were characterized by Gelfand and
Smirnov~\cite[pp. 194--197]{Gelfand-Smirnov} as the solutions of a certain
system of linear partial differential equations. In this section we derive
equations that are simpler, but which constitute a {\it sufficient condition\/}
for an integrand to be projective. In the two cases that interest us most, $k=1$
and $k = n-1$ they also constitute a necessary condition. It is the simplicity
of the equations that characterize projective $(n-1)$-integrands that will enable
us to prove our main result.

If $\Phi$ is a differential $k$-integrand, we denote the exterior
differential of the function $\Phi(\bx,\cdot)$ on $\Lambda^k \R^n$ by
$\delta \Phi(\bx,\cdot)$. For each $a$ in $\Lambda^k \R^n$ we consider
$\delta \Phi(\bx,a)$ as a vector in $\Lambda^k \R^{n*}$. For example,
if we set $k = 1$ and use coordinates $(x_1,\ldots,x_n;v_1,\ldots,v_n)$ on
$\R^n \times \Lambda^1 \R^n$, then
\begin{equation}
\delta \Phi = \sum_{i=1}^n \frac{\partial \Phi}{\partial v_i} \, dx_i .
\end{equation}
From the homogeneity of $\Phi$, it follows that
$\delta \Phi(\bx,a) \cdot a = \Phi(\bx,a)$ and that
$\delta \Phi(\bx,t a) = \delta \Phi(\bx,a)$ for $t > 0$.
When $\Phi$ is a differential $k$-form in the sense that
$\Phi(\bx,a) = \beta_{\bx} \cdot a$ for some differential $k$-form $\beta$
on $\R^n$, $\delta \Phi(\bx,a) = \beta_\bx$. In general, $\delta \Phi$ can
be seen as a differential $k$-form on $\R^n$ depending a parameter
$a$. It is then possible to compute its exterior differential in the $\bx$
variable $d_\bx \delta \Phi$. In the case $k = 1$ this would be
\begin{equation}\label{pre-Hamel}
d_\bx \delta \Phi = \sum_{1\leq i,j \leq n}
\frac{\partial^2 \Phi}{\partial x_i \partial v_j} \, dx_i \wedge dx_j .
\end{equation}

\begin{theorem}\label{generalized-Hamel}
A sufficient condition for a differential $k$-integrand $\Phi$ on $\R^n$ to be
projective is that $d_\bx \delta \Phi = \bzero$. Moreover, for $k = 1$ and
$n-1$ this condition is also necessary.
\end{theorem}

For example, from~(\ref{pre-Hamel}) we have that {\it a
differential $1$-integrand $\Phi$ on $\R^n$ is projective if and only if
it satisfies Hamel's equations\/} (\cite{Hamel:1903})
\begin{equation}\label{Hamel's-equations}
\frac{\partial^2 \Phi}{\partial x_i \partial v_j} =
\frac{\partial^2 \Phi}{\partial x_j \partial v_i}
\quad (1 \leq i, j \leq n).
\end{equation}

\begin{definition}
If $\Phi$ is a differential $k$-integrand in $\R^n$, its Hilbert-Lepage $k$-form is
the differential $k$-form $\alpha$ on $\R^n \times \Lambda^k \R^n$ defined by
$$
\alpha_{(\bx,a)}\left((\dot{\bx}_1,\dot{a}_1),\ldots,(\dot{\bx}_k,\dot{a}_k)\right)
=
\delta \Phi (\bx,a) \cdot (\kvector{\dot{\bx}}{k}) .
$$
\end{definition}

If we use coordinates $a_{i_1 \cdots i_k}$ in $\Lambda^{k} \R^n$ taken with
respect to its standard basis $\be_{i_1} \wedge \cdots \wedge \be_{i_k}$
$(1 \leq i_1 < \cdots < i_k \leq n)$,
the Hilbert-Lepage form can be written as
$$
\sum_{i_1 \cdots i_k} \frac{\partial \Phi}{\partial a_{i_1 \cdots i_k}} \,
dx_{i_1} \wedge \cdots \wedge dx_{i_k} .
$$

Let $N$ be  an oriented $k$-dimensional submanifold of $\R^n$. Given a nonzero
section $v$ of the line bundle $\Lambda^k(TN) \subset N \times \Lambda^k \R^n$
that represents the orientation of $N$, we can lift $N$ to
$\R^n \times \Lambda^k \R^n$ by defining
$$
\tilde{N} = \{(\bx,a) \in \R^n \times \Lambda^k \R^n : \bx \in N, a = v(\bx) \} .
$$
The key property of the Hilbert-Lepage form is that for any such lift,
$$
\int_{N} \Phi = \int_{\tilde{N}} \alpha .
$$
This identity allows us to derive the Euler-Lagrange equations in a very elegant
form:

\begin{proposition}\label{Euler-Lagrange}
A $k$-dimensional submanifold $N \subset \R^n$ is an extremal for the variational
problem posed by the differential $k$-integrand $\Phi$ if and only if for any
$k$-vector $b$ tangent to $\tilde{N}$ at a point $(\bx, a)$, we have that
$d\alpha_{(\bx,a)} \rfloor b = 0$.
\end{proposition}

\noindent {\it Sketch of the proof.}
Let $N_t$ be a variation of $N$ such that all the submanifolds $N_t$ agree
with $N$ outside a (small) compact set. If we lift $N_t$ to a variation
$\tilde{N}_t$ of $\tilde{N}$, we have
$$
\frac{d}{dt} \int_{N_t} \Phi \,{\Big|}_{t=0}  =
\frac{d}{dt} \int_{\tilde{N}_t} \alpha \, {\Big|}_{t=0} =
\int_{\tilde{N}} L_X \alpha ,
$$
where the Lie derivative $L_X$ is taken with respect to the vector field
$X$ along $N$ that determines the infinitesimal variation. Cartan's
formula and the fact that the variation $\tilde{N}_t$ is of compact support
imply that
$$
\int_{\tilde{N}} L_X \alpha = \int_{\tilde{N}} d\alpha \rfloor X +
\int_{\tilde{N}} d(\alpha \rfloor X) = \int_{\tilde{N}} d\alpha \rfloor X .
$$
Since the infinitesimal variation $X$ is an arbitrary vector field along
$\tilde{N}$ with compact support, it follows that
$d\alpha_{(\bx,a)} \rfloor b = 0$ for any $k$-vector $b$ tangent to $\tilde{N}$
at an arbitrary point $(\bx,a)$.
\qed \medskip

\noindent {\it Proof of Theorem~\ref{generalized-Hamel}.}
If $N$ is a $k$-dimensional affine subspace of $\R^n$, we can take its
lift $\tilde{N}$ to be the set of pairs $(\bx,a)$, where $\bx$ ranges over
all points on $N$ and $a$ is a {\it fixed\/} $k$-vector. Any $k$-vector $b$
tangent to $\tilde{N}$ at $(\bx,a)$ is of the form $(\lambda a,0)$,
where $\lambda$ is any real number. It follows easily from the definition
and the homogeneity of the Hilbert-Lepage form that
$$
d \alpha_{(\bx,a)} \rfloor (\lambda a, 0) =
d_{\bx} \delta \Phi(\bx,a) \rfloor \lambda a .
$$
Therefore, by Proposition~\ref{Euler-Lagrange}, if
$d_{\bx} \delta \Phi(\bx,a) = 0$ all $k$-planes are extremal for the variational
problem posed by the differential $k$-integrand $\Phi$.

The necessity of the equation $d_{\bx} \delta \Phi(\bx,a) = 0$ in the case
$k=1$ (i.e., Hamel's equations~\ref{Hamel's-equations}) is well known
(see Theorem 5.5.2 in \cite{Tabachnikov-Ovsienko} for a nice quick proof).
For $k = n-1$ the proof is even easier if we introduce coordinates
$(a_1,\ldots,a_n)$ in $\Lambda^{n-1} \R^n$ taken with respect
to the basis formed by the $(n-1)$-vectors
$$
(-1)^{i-1} \be_1 \wedge \cdots \wedge \be_{i-1} \wedge
\widehat{\be_{i}} \wedge \be_{i+1} \wedge \cdots \wedge \be_n
\quad (1 \leq i \leq n) .
$$
In these coordinates,
$$
d_{\bx} \delta \Phi =
\left( \sum_{i=1}^n
       \frac{\partial^2 \Phi}{\partial x_i \partial a_i}(\bx,a) \right)
\, dx_1 \wedge \cdots \wedge dx_n  ,
$$
and the condition $d_{\bx} \delta \Phi (\bx, a) \rfloor a = 0$  translates to
$$
\left( \sum_{i=1}^n
       \frac{\partial^2 \Phi}{\partial x_i \partial a_i}(\bx,a) \right) a_i = 0
\quad \mbox{for} \ 1 \leq i \leq n .
$$
This clearly implies that $d_{\bx} \delta \Phi = 0$ if and only if
$d_{\bx} \delta \Phi (\bx, a) \rfloor a = 0$ for all $(n-1)$-vectors $a$.
\qed \medskip

The last part of the preceding proof presents us with a very simple criterion
for determining whether a differential $(n-1)$-integrand is projective:

\begin{corollary}\label{Berck}
A differential $(n-1)$-integrand $\Phi$ on $\R^n$ is projective if and only if
it satisfies the equation
\begin{equation}\label{Berck's-equation}
\sum_{i=1}^n \frac{\partial^2 \Phi}{\partial x_i \partial a_i} = 0 .
\end{equation}
\end{corollary}

We're finally ready to prove our main result.

\begin{theorem}\label{Main-theorem}
For any value of the parameter $\lambda$, all geodesics of the Finsler metric
$$
\varphi_\lambda (\bx,\bv) =
\frac{(1 + \lambda^2 \nmx{\bx}^2)\nmx{\bv}^2 + \lambda^2 \innerxy{\bx}{\bv}^2}
{\nmx{\bv}}
$$
on $\R^3$ are straight lines. However, the only value of $\lambda$ for which
all planes are extremals of the Hausdorff $2$-area functional of $\varphi_\lambda$
is $\lambda = 0$.
\end{theorem}

\Proof
The proof of the first part of the theorem is reduced to showing that
$\varphi_\lambda$ satisfies Hamel's equations~\ref{Hamel's-equations} for
any value of $\lambda$. This is simplified by the linearity of the equations
and the decomposition of $\varphi_{\lambda}(\bx,\bv)$ into
$$
\nmx{\bv} + \lambda^2 \psi(\bx,\bv) =
\nmx{\bv} + \lambda^2
  \frac{\nmx{\bx}^2\nmx{\bv}^2 + \innerxy{\bx}{\bv}^2}{\nmx{\bv}} .
$$
It is easy to see that both $\psi$ and $\bv \mapsto \nmx{\bv}$ satisfy
Hamel's equations..

For the second part, we have to show that the Hausdorff $2$-area integrand
of $\varphi_\lambda$ (Theorem~\ref{area-integrand})
$$
\phi_{\lambda}(\bx,\ba) =
\frac{2(1 + \lambda^2\nmx{\bx}^2)^{3/2}
  \left((1 + 2\lambda^2\nmx{\bx}^2)\nmx{\ba}^2 -
                         \lambda^2 \innerxy{\bx}{\ba}^2 \right)^{3/2}}
 {(2 + 3 \lambda^2\nmx{\bx}^2)\nmx{\ba}^2 -  \lambda^2 \innerxy{\bx}{\ba}^2}
$$
does {\bf not} satisfy Equation~\ref{Berck's-equation} unless $\lambda$ is
zero. This is a nasty computation, but one that can easily be carried out
by machine. For example, using Maple we find that
$$
\frac{\partial^2 \phi_\lambda}{\partial x_1 \partial a_1} +
\frac{\partial^2 \phi_\lambda}{\partial x_2 \partial a_2} +
\frac{\partial^2 \phi_\lambda}{\partial x_3 \partial a_3}
$$
evaluated at the point $(t,t,t;1,1,0)$ in $\R^3 \times \Lambda^2\R^3$ equals
$$
-24 \frac{(1 + 3\lambda^2 t^2)^{3/2} \lambda^8 t^7}
       {(2 + 7\lambda^2 t^2)^3 \sqrt{2 + 8 \lambda^2 t^2}}.
$$
This is equal to zero for all values of $t$ if and only if $\lambda= 0$. In
this case, the integrand $\phi_\lambda$ is just the Euclidean area integrand
in $\R^3$ and all planes are minimal.
\qed \medskip

As was mentioned in the introduction, this theorem immediately implies that
if the Hausdorff measure is used as the notion of volume on Finsler spaces,
totally geodesic submanifolds are not necessarily extremal for the
area functional. By contrast, Berck has proved in \cite{Berck:2004} that
{\it if the Holmes-Thompson volume is used, then totally geodesic submanifolds
are extremal.}

The non-minimality of (totally geodesic) planes for the Finsler metric
$\varphi_\lambda$ for $\lambda \neq 0$ also shows that the Holmes-Thompson
volume cannot be replaced by the Hausdorff measure in the following
filling theorem of Ivanov.

\begin{theorem}[Ivanov, \cite{Ivanov:2002}]
Let $\varphi$ be a Finsler metric on the closed two-dimensional disc $D$
such that every two points on $D$ are joined by a unique geodesic.
If $\psi$ is another Finsler metric on $D$ such that the distance
induced by $\psi$ on the boundary $\partial D$ of $D$ is greater than
or equal to the distance induced by $\varphi$ on $\partial D$, then
the Holmes-Thompson volume of $(D,\varphi)$ does not exceed that of
$(D,\psi)$.
\end{theorem}

In fact, Ivanov's theorem has the following corollary:

\begin{corollary}
A totally geodesic two-dimensional submanifold $N$ of a Finsler manifold
$(M,\varphi)$ is minimal with respect to the Holmes-Thompson area integrand.
\end{corollary}

The way in which {\it minimal\/} differs from {\it extremal\/} is made clear
in the proof.

\Proof
Let $\bx$ be a point in $N$ and let $U_\bx$ be a small convex neighborhood of
$\bx$ in $M$. Such neighborhoods exist around every point by a well-known
theorem of J.~H.~C. Whitehead (\cite{Whitehead:1932}). The intersection of
$U_\bx$ and $N$ is a disc $D_\bx$ where every two points can be joined by a
unique geodesic. If $N'$ is another two-dimensional submanifold that coincides
with $N$ outside $U_\bx$ and such that the intersection $D'_\bx = U_\bx \cap N'$ is
a disc, then the (Holmes-Thompson) area of $N$ is less than or
equal to the area of $N'$.

To prove this, let $D$ be the two-dimensional disc and let $f$ and $f'$ be two
embeddings of $D$ into $M$ such that $f(D) = D_\bx$ and $f'(D) = D'_\bx$. The
Finsler metrics $f^* \varphi$ and $f'^{*}\varphi$ on $D$ satisfy the
hypotheses of Ivanov's theorem and, therefore, the area of $D_\bx$ does not
exceed the area of $D'_\bx$.
\qed \medskip

In the next section we prove a similar result for totally geodesic
hypersurfaces.

\section{Crofton formula for Finsler spaces}\label{Crofton-section}

Let $M$ be a Finsler manifold such that its space of oriented
geodesics is a manifold $\geod{M}$. Let $S^{*}M$ denote its unit
co-sphere bundle  and consider the canonical projection
$\pi : S^{*}M \rightarrow \geod{M}$ that sends a given unit covector
to the geodesic that has this covector as initial condition.

\begin{defiprop}[\cite{Arnold:1989} and \cite{Besse}]
Let $M$ be a Finsler manifold with manifold of geodesics $\geod{M}$
and let
$$
\begin{CD}
S^*M   @>i>>  T^*M \\
@V{\pi}VV          \\
\geod{M}
\end{CD}
$$
be the canonical projection onto $\geod{M}$ and the canonical inclusion
into $T^{*}M$. If $\omega_{\ssM}$ is the standard symplectic form on
$T^{*}M$, then there is a unique symplectic form $\omega$ on $\geod{M}$
which satisfies the equation $\pi^{*}\omega = i^{*}\omega_{\ssM}$.
\end{defiprop}

Among the many examples of Finsler manifolds whose space of geodesics is
a manifold we find convex neighborhoods in Finsler spaces,
Hadamard Riemannian manifolds (for both of these examples see
\cite{Ferrand:1997}), Zoll metrics (\cite{Besse}), Hilbert geometries,
and, more generally, projective Finsler metrics (\cite{Alvarez:symplectic}
and \cite{Alvarez-Fernandes:1998}). For all these manifolds we have the following
integral-geometric formula for the Holmes-Thompson area of hypersurfaces.

\begin{theorem}\label{crofton-hypersurface}
Let $M$ be an $n$-dimensional Finsler manifold with manifold of
geodesics $\geod{M}$. If $\omega$ is the symplectic form on $\geod{M}$
and  $N \subset M$ is any immersed hypersurface, then
$$
\vol_{n-1}(N) = \frac{1}{2\epsilon_{n-1}(n-1)!} \, \int_{\gamma \in \geod{M}}
\#(\gamma \cap N) |\omega^{n-1}| ,
$$
where $\epsilon_{n-1}$ is the volume of the Euclidean unit ball of dimension
$n-1$.
\end{theorem}

This theorem is probably folklore, but after its removal from the final
version of \cite{Alvarez:symplectic} a proof cannot be found in the
literature. The theorem follows from two basic lemmas, one based on the
coarea formula and the other based on symplectic reduction.

\begin{lemma}\label{coarea-lemma}
If $S^{*}M|_{\ssN}$ denotes the unit co-sphere
bundle of $M$ restricted to $N$, then
$$
\int_{S^{*}M|_{\ssN}} |\omega_{\ssM}^{n-1}|
= \int_{\gamma \in \geod{M}} \#(\gamma \cap N) |\omega^{n-1}|.
$$
\end{lemma}

\Proof
This is a simple application of the coarea formula:

{\it Let $X$ and $Y$ be smooth manifolds of the same dimension, let
$f : X \rightarrow Y$ be a smooth map, and let $\Omega$ be a
top-order differential form on $Y$. If for every regular value of
$f$ the number of preimages is finite, then}
\begin{equation}\label{coarea-formula}
\int_{X} f^{*} |\Omega| = \int_{y \in Y} \#(f^{-1}(y))\, |\Omega|.
\end{equation}

To prove the lemma, we just  have to set $X$ equal to $S^{*}M|_{\ssN}$,
$Y$ equal to $\geod{M}$, $f$ equal to the restriction of the projection
$\pi$ to $S^{*}M|_{\ssN}$, and use that, by
definition, $\pi^{*} \omega = i^{*}\omega_{\ssM}$.
\qed

\begin{lemma}\label{symplectic-lemma}
Let $\rho : S^{*}M|_{\ssN} \rightarrow T^{*}N$ be the map that takes a
covector $\xi_{n}$ in $S^{*}M|_{\ssN}$ and sends it to the covector
$\xi_{n}|_{T_{n}N}$ in $T^{*}N$. If $j : S^{*}M|_{\ssN} \rightarrow T^{*}M$
is the canonical inclusion, then $\rho^{*}\omega_{\ssN} = j^{*} \omega_{\ssM}$.
\end{lemma}

\Proof
The submanifold $T^{*}M|_{\ssN} \subset T^{*}M$ is a coisotropic submanifold
and the leaves of its characteristic foliation are the conormals of
$N$. The reduced phase space is identified with $T^{*}N$ via the
projection $\xi_{n} \mapsto \xi_{n}|_{T_{n}N}$ and it is easy to see that
the reduced symplectic form coincides with the standard form on
$T^{*}N$. From this follows that $\rho^{*}\omega_{\ssN} = j^{*} \omega_{\ssM}$.
\qed \medskip

Notice that the image of $ S^{*}M|_{\ssN}$ under $\rho$ equals the
(closed) unit co-disc bundle $B^{*}N$ of the submanifold $N$ with its induced
Finsler metric. Moreover, every point in the interior of $B^{*}N$ is a regular
value and has exactly two preimages.

\noindent {\it Proof of Theorem~\ref{crofton-hypersurface}.}
Using the previous remark and applying the coarea formula, we have that
$$
\int_{S^{*}M|_{\ssN}} \rho^* |\omega_{\ssN}^{n-1}| =
2 \int_{B^{*}N} |\omega_{\ssN}^{n-1}| = 2\epsilon_{n-1}(n-1)! \, \vol_{n-1}(N) .
$$
The last equality is simply the definition of the Holmes-Thompson volume of
$N$.

Applying Lemmas~\ref{symplectic-lemma} and~\ref{coarea-lemma}, we obtain
\begin{eqnarray*}
\vol_{n-1}(N) &=&
\frac{1}{2\epsilon_{n-1}(n-1)!} \,
\int_{S^{*}M|_{\ssN}} \rho^* |\omega_{\ssN}^{n-1}| \\
&=& \frac{1}{2\epsilon_{n-1}(n-1)!} \,
    \int_{S^{*}M|_{\ssN}} |\omega_{\ssM}^{n-1}| \\
&=& \frac{1}{2\epsilon_{n-1}(n-1)!} \,
    \int_{\gamma \in \geod{M}} \#(\gamma \cap N) |\omega^{n-1}|.
\end{eqnarray*}
\qed \medskip

The importance of this theorem is underlined by the following application:

\begin{theorem}
A totally geodesic hypersurface $N$ of a Finsler manifold $(M,\varphi)$
is minimal with respect to the Holmes-Thompson area integrand.
\end{theorem}

\Proof
Let $\bx$ be a point in $N$ and let $U_\bx$ be a convex neighborhood of $\bx$
in $M$. The space of geodesics of $U_\bx$ with the restricted Finsler metric is
easily seen to be a manifold. In fact, since any geodesic passing through
$U_\bx$ intersects its boundary in two distinct points, the space of (oriented)
geodesics of $U_\bx$ can be identified with
$\partial U_\bx \times \partial U_\bx$ minus the diagonal. If $N'$ is
another hypersurface that coincides with $N$ outside of $U_\bx$, then the
(Holmes-Thompson) area of $N$ is less than or equal to the area of $N'$.

The proof depends on two simple remarks on the behavior of geodesics in
$U_\bx$: (a) every geodesic that intersects $N \cap U_\bx$ is either contained
in it or intersects it in only one point; (b) any geodesic that intersects
$N \cap U_\bx$ intersects $N' \cap U_\bx$.

The first remark follows from the convexity of $U_\bx$ and the fact that $N$
is totally geodesic. Indeed, the existence of a geodesic $\gamma$ not lying on
$N \cap U_\bx$ and intersecting it at two points would imply the existence of
at least two geodesics in $U_\bx$ joining these points:
the geodesic $\gamma$ and some geodesic on $N \cap U_\bx$. The second remark
is topological: the union of $N \cap U_\bx$ and $N' \cap U_\bx$ is a (non-smooth)
hypersurface that divides its complement in $U_\bx$ into an inside and an outside.
Since a geodesic in $U_\bx$ can cross $N \cap U_\bx$ at most once; if it crosses
at all, it must also cross $N' \cap U_\bx$. In other words, except for a set of
measure zero (those geodesics lying on $N \cap U_\bx$) the geodesics $\gamma$ in
$\geod{U_\bx}$ satisfy
$\#\left((N'\cap U_\bx) \cap \gamma)\right) \geq
\#\left((N\cap U_\bx) \cap \gamma)\right)$. Therefore, by
Theorem~\ref{crofton-hypersurface}, the Holmes-Thompson area of $N \cap U_\bx$
is less than or equal to that of $N' \cap U_\bx$.
\qed \medskip

As the previous theorem shows, there is a very close relation between minimality
of totally geodesic hypersurfaces and Crofton formulas. Another instance of this
phenomenon is the following characterization by Gelfand, Smirnov, \'Alvarez,
and Fernandes (see Theorem~4 of \cite{Gelfand-Smirnov} and Theorem~6.1
of~\cite{Alvarez-Fernandes:2000} for a rectification and proof) of {\it even\/}
projective differential $(n-1)$-integrands. We use the term {\it even\/} to
denote integrands $\Phi$ that satisfy $\Phi(\bx,-a) = \Phi(\bx,a)$.

\begin{theorem}
An even differential $(n-1)$-integrand $\Phi$ on $\R^n$ is projective if and only
if there exists a smooth (possibly signed) measure $\mu$ on the space $\geod{\R^n}$
of oriented lines on $\R^n$ such that if $N$ is any hypersurface,
$$
\int_N \Phi = \int_{\ell \in \geod{\R^n}} \#(\ell \cap N) \, d\mu .
$$
\end{theorem}

In other words, the existence of a Crofton formula involving the space of lines
in $\R^n$ is a necessary and sufficient condition for hyperplanes to be extremal.
This remark, together with Theorem~\ref{Main-theorem}, gives us the following
result:

\begin{corollary}
The only value of the parameter $\lambda$ for which the Hausdorff $2$-area
integrand of the projective Finsler metric
$$
\varphi_\lambda (\bx,\bv) =
\frac{(1 + \lambda^2 \nmx{\bx}^2)\nmx{\bv}^2 + \lambda^2 \innerxy{\bx}{\bv}^2}
{\nmx{\bv}}
$$
on $\R^3$ admits a Crofton formula is $\lambda = 0$.
\end{corollary}

By restricting the metrics $\varphi_\lambda$ to a compact ball $B$ and letting
$\lambda$ tend to zero, we obtain projective metrics that are arbitrarily
close to the Euclidean metric on $B$ and for which, if the Hausdorff measure is
employed, there is no Crofton formula for hypersurfaces.

While our results suggest that the Hausdorff measure is generally unsuitable
for the study of Finsler spaces, its study in this context may lead to
interesting rigidity theorems. We have in mind results such as those of
Colbois, Vernicos, and Verovic who prove in \cite{Colbois-etal:2003} that
{\it if the Hausdorff $2$-dimensional area of all ideal triangles in a
Hilbert geometry is a fixed constant, then it is the classical hyperbolic
geometry.} Their result and those of the present work suggest the following
question:

\begin{problem}
Does the fact that hyperplanes minimize the Hausdorff area functional characterize
hyperbolic geometry among all the Hilbert geometries?
\end{problem}

\medskip
\centerline{\sc Acknowledgements}
\medskip
The first author is happy to acknowledge the hospitality and great working
conditions provided by the Universit\'e Catholique de Louvain, the
Centre de Recerca Matem\`atica, and the UNICAMP during the different stages
of the preparation of this article.

\bibliography{../../paperbib}
\bibliographystyle{amsplain}

\end{document}